\documentclass[12pt]{article}

\usepackage{amsmath,amssymb,amsfonts,theorem,makeidx,latexsym,epsfig,,subfigure}

\newtheorem{defn}{Definition}[section]

\newtheorem{lemma}[defn]{Lemma}

{\theorembodyfont{\rmfamily}

\newtheorem{ex}[defn]{Example}}

\newtheorem{thm}[defn]{Theorem}

\newtheorem{prop}[defn]{Proposition}

\newtheorem{cor}[defn]{Corollary}

\numberwithin{equation}{section}

\newcommand{\h}{{\cal H}}

\newcommand{\ltr}{ L^2(\mathbb R) }

\newcommand{\ltn}{{\ell}^2(\mathbb N)}

\newcommand{\ltz}{{\ell}^2(\mathbb Z)}

\newcommand{\si}{S^{-1}}

\newcommand{\mn}{\mathbb N}

\newcommand{\mr}{\mathbb R}

\newcommand{\mz}{\mathbb Z}

\newcommand{\mc}{\mathbb C}

\newcommand{\mts}{ \{E_{mb}T_{na}g \}_{m,n \in \mz}}

\def\bp{{\noindent\bf Proof. \ }}

\def\ep{\hfill$\square$\par\bigskip}

\def\bqs{\begin{equation}}

\def\eqs{\tag*{$\square$}\end{equation}\par\bigskip}

\def\la{\langle}

\def\ra{\rangle}

\def\ga{\gamma}

\def\ftk{\{f_k\}_{k=1}^\infty}

\def\ctk{\{c_k\}_{k=1}^\infty}

\def\ctz{\{c_k\}_{k\in \mz}}

\def\gtk{\{g_k \}_{k=1}^\infty}

\def\etk{\{e_k\}_{k=1}^\infty}

\def\suk{\sum_{k=1}^\infty}

\def\sukz{\sum_{k\in \mz}}

\def\nl{\left|\left|}

\def\nr{\right|\right|}

\def\span{\overline{\text{span}}}

\def\Span{\text{span}}

\def\supp{\text{supp}}

\def\bop{\begin{op}\rm}

\def\eop{\end{op}}

\def\bee{\begin{eqnarray}}

\def\ene{\end{eqnarray}}

\def\bes{\begin{eqnarray*}}

\def\ens{\end{eqnarray*}}

\def\bei{\begin{itemize}}

\def\eni{\end{itemize}}

\def\bt{\begin{thm}}

\def\et{\end{thm}}

\def\bc{\begin{cor}}

\def\ec{\end{cor}}

\def\bpr{\begin{prop}}

\def\epr{\end{prop}}

\def\bl{\begin{lemma}}

\def\el{\end{lemma}}

\def\bd{\begin{defn}}

\def\ed{\end{defn}}

\def\bex{\begin{ex}}

\def\enx{\end{ex}}

\def\bfi{\begin{fig}}

\def\efi{\end{fig}}

\def\inr{\int_{-\infty}^\infty}

\def\sukz{\sum_{k\in \mz}}

\def\wt{\widetilde{T}}

\def\ftz{\{f_k\}_{k\in \mz}}
\def\gtz{\{g_k\}_{k\in \mz}}
\newcommand{\Tnz}{\{T^kf_0\}_{k\in\mz}}
\newcommand{\Vnz}{\{V^kg_0\}_{k\in\mz}}
\newcommand*{\numbersys}[1]{\ensuremath{\mathbb{#1}}}

\newcommand*{\R}{\numbersys{R}}

\newcommand*{\Z}{\numbersys{Z}}

\title{Operator representations of frames: boundedness, duality, and stability. }
\date{\today}

\author{Ole Christensen, Marzieh Hasannasab }

\begin{document}
	
	\maketitle

	\begin{abstract}
The purpose of the paper is to analyze frames $\ftz$ having the form $\{T^kf_0\}_{k\in\mz}$ for some linear operator $T: \mbox{span} \ftz \to \mbox{span}\ftz$.  A key result
characterizes boundedness of the operator $T$ in terms of shift-invariance of a
certain sequence space. One of the consequences is
 a characterization of the case where the representation
$\ftz=\{T^kf_0\}_{k\in\mz}$
 can be achieved for an operator
$T$ that has an extension to a bounded bijective operator $\widetilde{T}: \h \to \h.$ In this
case we also characterize all the dual frames that are representable in terms
of iterations of an operator $V;$ in particular we prove that the only possible operator
is $V=(\widetilde{T}^*)^{-1}.$ Finally, we consider stability of the representation $\{T^kf_0\}_{k\in\mz};$
rather surprisingly, it turns out that the possibility to represent a frame on
this form is sensitive towards some of the
classical perturbation conditions in frame theory. Various ways of avoiding this problem
will be discussed.
Throughout the paper the results will be connected with the
operators and function systems appearing in applied harmonic analysis, as well as with general group
representations.
%The developed theory covers not only the mentioned cases of
%operators on $L^2(\R)$ but also natural generalizations to LCA groups and other classes of operators.
	\end{abstract}
	\section[Results]{Introduction}
In this paper we consider frames
$\ftz$ in a Hilbert space $\h$ arising via iterated action of a linear operator $T: \Span \ftz \to
\Span \ftz,$ i.e., on the form
\bee \label{61019a} \ftz=\Tnz.\ene We say that the frame $\ftz$ is {\it represented} via the operator $T.$
The motivation to consider frames of this form comes from several directions:
\bei
\item The (Fourier) orthonormal basis $\{e^{2\pi ikx}\}_{k\in \mz}$ for $L^2(0,1)$ has the form
\eqref{61019a}, where $(Tf)(x)=e^{2\pi ix} f(x)$ and $f= \chi_{[0,1]}.$
\item Single-generated shift-invariant systems (Example \ref{motivation example}) and
Gabor systems (Example \ref{60829b}) have the form \eqref{61019a}.
\item A group representation
acting on a cyclic (sub)group indexed by $\mz$ (Example \ref{61006a}) leads to a system of vectors on the from \eqref{61019a}.
\eni

The idea of representing frames
on the form \eqref{61019a} is also closely related with dynamical sampling, see, e.g., \cite{A1}. However,
the indexing of a frame in the context of dynamical sampling  is different from the one used in
\eqref{61019a}, and we will show that a re-indexing might change the properties of the operator $T$ drastically.

In Section \ref{61019c} we first classify the frames having the form
$\ftz=\Tnz,$ where $T$ is
a linear (not necessarily bounded) operator defined on $\mbox{span}\ftz.$ One of the main results characterizes the frames
that can be represented in terms of a bounded operator $T,$ in terms of shift-invariance of
 a certain subspace of $\ell^2(\mz).$ Various consequences of this result are derived, e.g.,  that if an overcomplete frame with finite excess
has a representation of the form \eqref{61019a}, then $T$ is necessarily
unbounded.

Section \ref{61003b} deals with the properties of the dual frames associated with
a frame on the form $\ftz=\Tnz.$ For the important case where $T$ is bounded and bijective,
we characterize the dual frames that can be represented in terms of a bounded operator $V;$ in particular, we
show that the only possibility of the representing operator
for the dual frame is $V=(T^*)^{-1}.$

In
Section \ref{61003c} we consider stability of a representation
\eqref{61019a} under various perturbation conditions. Rather surprisingly, it turns out that
a representation of such a type is unstable under the classical perturbation conditions
in frame theory, e.g., the Paley-Wiener type conditions \cite{C13,OC}. That is, a perturbation of
a frame $\ftz=\Tnz$ might not be representable in terms of an operator; or, if the operator $T$
is bounded, a perturbation might turn the frame $\ftz$  into  a frame that is only representable in terms of an unbounded operator. We prove, however, that under certain restrictions on the
perturbation condition,  stability and boundedness is preserved. Finally, for frames $\ftz$
that are norm-bounded below we prove that the type of  perturbation condition
that is used most frequently in the literature leads to frames $\gtz$ that
can be represented via iterations of a {\it finite} collection of bounded operators.

The paper closes with an appendix, containing some operator theoretic considerations.  We show, e.g.,
that the chosen indexing is important for the properties of the operator representing a given
frame.

In the rest of the introduction, we will collect some definitions and standard results from frame
theory. A sequence $\ftk$
in a Hilbert space $\h$ is a {\it frame}
for  $\h$ if there exist constants $A,B>0$ such that
$ A \, ||f||^2 \le \suk | \la f, f_k\ra|^2 \le B \, ||f||^2, \, \forall f\in \h.$
The sequence $\ftk$ is a Bessel sequence if at least the upper frame condition holds.  Also, $\ftk$ is called a {\it Riesz sequence} if
there exist constants $A,B>0$ such that
$A \sum |c_k|^2 \le \nl \sum c_k f_k \nr^2 \le B\sum |c_k|^2$
for all finite scalar sequences $c_1, \dots, c_N, \, N\in \mn.$

If $\ftz$ is a Bessel sequence, the {\it synthesis operator} is defined by
\bee \label{60811a} U: \ltz \to \h, \, U\ctz := \sukz c_k f_k;\ene
it is well known that $U$ is well-defined and bounded. A central role will be played by the kernel
of the operator $U,$ i.e., the subset of $\ltz$ given by
\bee \label{60811f} N_U=\left\{\ctz\in\ltz~\bigg|~\sum_{k\in\Z}c_kf_k=0\right\}.\ene
The {\it excess} of a frame  is the number of elements that can be removed yet leaving a frame.
It is well-known that the excess equals $\mbox{dim} (N_U);$ see  \cite{Balan}.

Given a Bessel sequence $\ftz,$ the {\it frame operator} $S: \h \to \h$ is defined by $S:=UU^*.$
For a frame $\ftz,$ the frame operator is invertible and
$f= \sum_{k\in \mz} \la f, \si f_k\ra f_k, \, \forall f\in \h.$
The sequence $\{\si f_k\}_{k\in \mz}$ is also a frame; it is called the {\it canonical dual frame.}

 For a frame $\ftz$ that is not a Riesz basis, it is known that there exists infinitely
many {\it dual frames}, i.e., frames $\gtz$ such that
\bes f= \sum_{k\in \mz} \la f, g_k\ra f_k, \, \forall f\in \h.\ens
The class of dual frames have been characterized by Li \cite{Li}.

Throughout the paper we will illustrate the results with applications to
frames appearing in applied harmonic analysis, e.g., shift-invariant
systems and Gabor systems.  First, for $a\in \mr,$ define the {\it translation operator}
$T_a$ acting on $\ltr$ by $T_af(x):= f(x-a)$ and the {\it modulation operator}
$E_a$ by $E_af(x):=e^{2\pi i ax} f(x).$  Both operators are unitary.  Furthermore,
defining the Fourier transform of $f\in L^1(\mr)$ by $\widehat{f}(\ga)={\cal F}f(\ga)=
\inr f(x)e^{-2\pi i \ga x}dx$ and extend it in the standard way to a unitary operator
on $\ltr,$ we have ${\cal F}T_a=E_{-a}{\cal F}.$
The following example will inspire us throughout the paper.
\begin{ex}\label{motivation example}
Consider a function $\varphi\in L^2(\R)\setminus \{0\}.$ Then the {\it shift-invariant system} $ \{T_k\varphi\}_{k\in\mz} $ is linearly independent. Letting
$\Phi(\ga):= \sum_{k\in \mz}|\widehat{\varphi}(\ga + k)|^2,$ it was proved in
\cite{BL1} (or see Theorem 9.2.5 in \cite{CB}) that   $ \{T_k\varphi\}_{k\in\mz} $ is a frame sequence if and only if
there exist $A,B>0$ such that $A\le\Phi(\ga)\le B,$ a.e. $\ga\in [0,1]\setminus N,$
where $ N=:\{ \ga\in [0,1]~ \big| ~  \sum_{k\in \mz}|\widehat{\varphi}(\ga + k)|^2 =0    \}.$
Furthermore, the special case where $\{T_k\varphi\}_{k\in\mz}$ is a Riesz sequence corresponds
to the case where the set $N$ has measure zero.
Note that regardless of the frame properties of the shift-invariant system,
we can write it on the form $\{T_k\varphi\}_{k\in\mz}= \{(T_1)^k \varphi\}_{k\in\mz},$
i.e., as the iterated system arising by letting the powers of the bounded
operator $T_1$ act on the function $\varphi.$ \ep \enx
%We also note that the operator $T_1$ is unitary; thus, the example shows that a system
%$\ftz=\Tnz$ can be overcomplete even for a unitary operator $T.$
%\ep \end{ex}
More generally, iterated systems $\Tnz$ naturally shows up in the context of group representations.
This topic is well connected with frame theory; see, e.g., the paper \cite{BHP} and the
references therein.

\bex \label{61006a} Let $G$ denote a locally compact group, and $\pi$ a group representation of $G$ on a Hilbert space $\h;$ that
is $\pi$ is a mapping from $G$ into the space of bounded invertible  operators on $\h,$ satisfying that
$\pi(xy)= \pi(x) \pi (y)$ forall $x,y\in G.$
Now, fix some $x_0\in G.$ Considering the cyclic subgroup of $G$ generated by the
element $x_0,$ i.e., the set $\{x_k\}_{k\in \mz} =\{x_0^k\}_{k\in \mz}\subset G,$  the group representation
acting  on a fixed $f_0\in \h$ generates the family of vectors
$ \{  \pi (x_0^k)f_0 \}_{k\in \mz} =  \{ [ \pi(x_0)]^k f_0\}_{k\in \mz}.$
This system has  the structure \eqref{61019a} with $T=\pi(x_0).$
Note that Example \ref{motivation example}
is a special case of this; indeed, the left-regular representation of the group $\mr$ with the composition " $ +$ "
on $\ltr$ is precisely
$(\pi(x)f)(t)= f(t-x), \, f\in \ltr, \, t,x\in \mr.$ The general setting of group representations
covers this example and its discrete variant in $\ltz$ (whose frame
properties are analyzed  in \cite{Jan2}) in a unifying way.

Note that the structure of systems arising
from a group representation is very rigid: arbitrary small perturbations might destroy the
special structure, so it is important that such cases can still be handled within the
frame work of the more general systems \eqref{61019a}.
\ep \enx

\section{The representation $\ftz=\Tnz$} \label{61019c}
In this section we want to consider representation of a frame $\ftz$ on the form $\Tnz$ for some linear operator $T$ defined on an appropriate subspace of $\h$. The starting point must be a clarification of the exact meaning of $\Tnz$. For $k\geq 0$, this is clear. For $k=-1$ we will require that $T$ is invertible as a map from $\Span\ftz$ into itself. This guarantees that $T^{-1}f_0$ is well-defined, and hence also $T^kf_0=(T^{-1})^{-k}f_0$ is well-defined for $k=-2,-3,\cdots$. In the following result, we characterize the availability of the representation $\ftz=\Tnz$. The proof is a modification of the corresponding result for sequences indexed by $\mn$, so we only sketch it.
\bpr \label{60811d}
	 	Consider a  frame sequence $\ftz$ in a Hilbert space $\h$ which spans an infinite dimensional subspace. The following are equivalent:
	 	\begin{itemize}
	 		\item[(i)] $\ftz$ is linearly independent.
	 		\item[(ii)] The map $Tf_k:= f_{k+1}$ is well-defined, and extends to a linear and invertible operator  $T:\Span\ftz\rightarrow \Span\ftz$.
	 	\end{itemize}
	 	In the affirmative case,  $ \ftz=\Tnz. $
	 \epr
	 \bp  The proof that $(i)\Rightarrow(ii) $ is exact the same as for sequences indexed by $\mn$, see \cite{olemmaarzieh}. Now assume that $(ii)$ holds. It is easy to see that $f_k\neq 0$, for all $k\in\mz$. Now in order to reach a contradiction, assume that $\sum_{k=M}^{N}c_kf_k=0$. For some coefficient $c_k$, $k=M,\cdots,N$ not all of which are zero. We can choose $M,N\in\mz$ such that $c_M\neq 0$, $c_N\neq 0$. Then, the same proof as in \cite{olemmaarzieh}
 shows that the vector space $V:=\Span\{f_k\}_{k=M}^{N}$ is invariant under the action of $T$. Now, a similar calculation shows that $V$ is invariant under the action of $T^{-1}$. Thus
 $\Span\ftz=\Span\Tnz=\Span\{f_k\}_{k=M}^{N}=V$, which is a contradiction because $\Span\ftz$ is assumed to be infinite-dimensional. Thus $\ftz$ is linearly independent, as desired.
 \ep

If$\ftz$ is a frame sequence and the operator $T$ in Proposition \ref{60811d} is bounded, it has a unique extension to a bounded operator $\widetilde{T}:
\span \ftz \to \span \ftz,$ given by
\bes \widetilde{T} \sukz c_k f_k = \sukz c_k f_{k+1}, \, \ctz\in \ltz.\ens
We first state a necessary condition in order for a frame $\ftz$ to have a representation
on the form $\Tnz$ for a given  bounded operator $T.$
	\bpr\label{lowerboundofT}
Consider a frame on the form $\ftz=\Tnz$ for some bounded linear operator
$T: \Span \ftz \to \Span \ftz$. Then the following hold:
\bei
\item[(i)] $\|T\|\geq 1$.
\item[(ii)] If $T^{-1}:\Span\ftz\rightarrow \Span\ftz$ is bounded, then $\|T^{-1}\|\geq 1$.
\eni
 	\epr
	\bp  Let $A,B$ denote some frame bounds and fix  any $n\in\mn$. Using the frame inequalities for any $f\neq 0$, we have
	\bee\label{1019be} A\|f\|^2 &\leq& \sum_{k\in\mz}|\langle f,T^kf_0\rangle|^2 =\sum_{k\in\mz}|\langle (T^n)^*f,T^{k-n}f_0\rangle|^2\nonumber\\
	& =& \sum_{k\in\mz}|\langle (T^n)^*f,T^{k}f_0\rangle|^2
	  \leq   B\|(T^n)^*f\|^2 \leq B\|T\|^{2n} \| f \|^2.
	\ene
	Therefore $A\leq B\|T\|^{2n}$ for all $n\in \mn,$ which implies that $\|T\|\geq 1$. The result in (ii) follows by  replacing $T$ by $T^{-1}$
and noticing that these two operators represent the same frame.\ep

Assuming that a frame $\ftz$ has a representation on the form $\Tnz,$  we will now characterize
boundedness of the operator $T$ in terms of the kernel of the synthesis operator, see
\eqref{60811f}; in particular, this leads to a characterization of the
case where the operator $T$  has an extension to a bounded bijective operator on $\h.$

For this purpose we need the analogue of the translation operator, acting on the sequence
space $\ltz.$ Define the {\it right-shift operator} on $\ltz$ by
\bee \label{60811b} {\cal T}:\ltz\rightarrow\ltz, \, {\cal T }\ctz=\{c_{k-1}\}_{k\in\mz}.\ene
Clearly ${\cal T}$ is a unitary operator on $\ltz$. We say that
a subspace $V\subseteq\ltz$ \textit{is invariant under right-shifts} (respectively, left-shifts) if ${\cal T}(V)\subseteq V$ (respectively, if ${\cal T}^{-1}(V)\subseteq V$).

			\begin{thm}\label{boundedness}
				Consider a  frame having the form $\ftz=\Tnz$
for some linear operator $T: \Span \ftz\to \Span \ftz,$ and let
  $A,B>0$ denote some frame bounds. Then the following hold:
  \bei
  \item[(i)] The operator $T$ is bounded if and only if the kernel $N_U$
  of the synthesis operator $U$ is invariant  under right-shifts; in the affirmative case,
   \bes 1 \le \| T\|\le \sqrt{BA^{-1}}.\ens
      \item[(ii)] The operator $T^{-1}:\Span\ftz\rightarrow\Span\ftz$ is bounded if and only if $N_U$ is invariant under left-shifts; in the affirmative case, \bes 1 \le \| T^{-1}\|\le \sqrt{BA^{-1}}.\ens
\item[(iii)] Assume that  $N_U$ is invariant under right and left-shifts.
         Then the operator $T$ has an extension to a bounded
         bijective operator $\widetilde{T}:\h \to \h.$
\eni \et

			\bp In order to prove (i), assume that $\{ f_k\}_{k\in\mz}= \{T^kf_0\}_{k\in\mz}$ for a bounded operator $T:\Span\ftz\rightarrow\Span \ftz.$  Then $T$ can be extended to a bounded operator $\widetilde{T} :\h \to\h$. For any $\ctz\in N_U$, we have
\bes U{\cal T}\ctz= \sum_{k=-\infty}^{\infty}c_{k-1}f_{k}= \sum_{k=-\infty}^{\infty}c_kf_{k+1}=\widetilde{T}\sum_{k=-\infty}^{\infty}c_kf_{k}=T0=0.\ens Therefore ${\cal T}(N_U)\subseteq N_U$, as claimed.
				
Conversely, assume that  $N_U$ is invariant  under right-shifts.
Assume that $f\in \Span \{ f_k\}_{k\in\mz},$ i.e., $f=\sum_{k=M}^{N}c_kf_k$ for some
$M,N\in \mz, c_k \in \mc.$  One can consider $\{c_k\}_{k=M}^N$ as a sequence $\{ c_k\}_{k\in\mz}$ in $\ltz$ where $c_k=0$ for $k>N$ and $k<M$. Thus we can write $\ctz=\{d_k\}_{k\in\Z}+\{r_k\}_{k\in\Z}$, where $\{d_k\}_{k\in\Z}\in N_U$ and $\{r_k\}_{k\in\Z}\in N_U^{\perp}$. Since $N_U$ is invariant under right-shifts, we have
				$\sum_{k\in\Z}d_kf_{k+1}=0$. Using the splitting of $\ctz$ and that $\{ f_k\}_{k\in\mz}$ is a Bessel sequence, we get that
				\bes \left\|Tf\right\|^2=\left\|\sum_{k=M}^{N}c_kf_{k+1}\right\|^2=
\left\|\sum_{k=-\infty}^{\infty}(d_k+r_k)f_{k+1}\right\|^2 & = &
\left\|\sum_{k=-\infty}^{\infty}r_k f_{k+1}\right\|^2 \\
				& \leq & B\sum_{k=-\infty}^{\infty}|r_k|^2. \ens
Recall (see Lemma 5.5.5 in \cite{CB}) that since $\ftz$  is a frame with lower bound
$A,$ we have
$A \sukz |c_k|^2 \le ||U\ctz||^2, \,  \forall \ctz \in N_U^\bot.$ It follows that
				\bes
				\left\|Tf\right\|^2
				\leq \frac{B}{A} \left\|\sum_{k\in\mz}r_kf_{k}\right\|^2  =  \frac{B}{A} \left\|\sum_{k\in\mz}c_kf_{k}\right\|^2
				=  \frac{B}{A}\left\|f\right\|^2,
				\ens
				%				Similar argument shows that $AB^{-1}\|f\|^2\leq \|\widetilde{T}(f)\|^2$.
				% Hence, \[\sqrt{AB^{-1}}\|f\|\leq\|\widetilde{T}f\|\leq \sqrt{A^{-1}B}\|f\|\qquad f\in\Span \ftz.\]
i.e., $T$ is bounded as desired. The above calculations also confirm the claimed upper
bound on the norm of $T.$ The lower bound in the estimate in (i) was proved in Proposition \ref{lowerboundofT}; this completes the proof of (i).
The result (ii) is a consequence of (i). Indeed, since $\ftz=\Tnz$, we can write $\{f_{-k}\}_{k\in\Z}=\{(T^{-1})^kf_0\}_{k\in\Z}$. Denoting
the synthesis operator for $\{f_{-k}\}_{k\in\Z}$ by $V,$
Theorem \ref{boundedness} shows that  $T^{-1}$ is bounded if and only if
the kernel $N_V$ is right-shifts invariant. It is easy to see that $\ctz\in N_U$ if and only if $\{c_{-k}\}_{k\in\Z}\in N_V$. Hence the left-shifts invariance  of $N_U$ is equivalent with the right-shift invariance of $N_V$.

For the proof of (iii),  if $N_U$ is invariant under right and left-shifts, then the operators $T,T^{-1}:\Span\ftz\rightarrow\Span\ftz$ are bounded. Hence they can be extended to
bounded operators $\widetilde{T}, \widetilde{T^{-1}}$ on $\h.$ Since
				\[T^{-1}Tf=TT^{-1}f=f, \quad f\in\Span\ftz,\]
				 it follows that $\widetilde{T}\widetilde{T^{-1}}= \widetilde{T^{-1}}\widetilde{T=}I,$
i.e., $\widetilde{T}$ is invertible on $\h$. \ep

Throughout the
paper it will be crucial to distinguish carefully between a bounded operator
$T: \Span \ftz\to \Span \ftz$ and its extension $\widetilde{T}: \h \to \h.$ Indeed, our setup
implies that $T$ is invertible,  but the extension to an operator on $\h$ might no longer
be injective (for the convenience of the interested reader we include such an example
in the Appendix).

Note that the biimplications in Theorem \ref{boundedness} uses the full strength of the frame
assumption.  Indeed, one can construct examples of sequences $\ftz= \Tnz$ satisfying only the upper
frame condition (resp. the lower frame condition), and such that $T$ is unbounded while the kernel $N_U$ is
invariant under right-shifts.

Let us demonstrate the power of Theorem \ref{boundedness} by some
consequences and examples; another application
will be given in Proposition \ref{61001b}.   Let us first consider the special case of a
Riesz sequence.

\bc \label{60829d} Any Riesz sequence $\ftz$ has a representation
$\Tnz$ for a bounded and bijective operator $T: \span \ftz \to \span \ftz.$ \ec

Corollary \ref{60829d} follows immediately from Theorem \ref{boundedness} and the
fact that the synthesis operator for a Riesz sequence  is injective. We therefore now
turn to the setting of an overcomplete frame.
			
\begin{cor} \label{60829a} Consider an overcomplete frame on the form $\ftz=\Tnz.$
If $T\in B(\h)$, then $ \mbox{dim}(N_U)=\infty.$
\end{cor}
\bp
If there is a nonzero element $c=\ctz$ in $N_U$, then by Theorem \ref{boundedness}, the boundedness of $T$ implies that ${\cal T}^j c\in N_U$ for all $j\in\mn$. We will now
show that the sequence $\{{\cal T}^j c\}_{j=1}^\infty$ is linearly independent; this
implies that  $N_U$ is  infinite-dimensional and concludes the proof. Now
consider the operator ${\cal F}:\ltz\mapsto L^2[0,1], \,  {\cal F}c=\sum_{k\in\mz}c_ke_k$, where $e_k(x)=e^{-2\pi i kx}$. The operator $\cal F$ is unitary, and $ {\cal F}{\cal T}c =e_{1} {\cal F}c$. Now assume that for some $N\in\mn$ and $d_1,d_2,\cdots,d_N\in\mc$,	we have $\sum_{j=1}^{N}d_j{\cal T}^jc=0$. Let $\varphi={\cal F}c$. Then
\[0={\cal F}(\sum_{j=1}^{N}d_j{\cal T}^j c)=\sum_{j=1}^{N}d_je_j{\cal F} c=(\sum_{j=1}^{N}d_je_j)\varphi.\]
This means that $  (\sum_{j=1}^{N}d_je_j(x))\varphi(x)=0 $, for $a. e.$ $ x\in[0,1] $. Since $\varphi\neq 0$, the support
of
$\varphi$   has positive measure. Thus we have $ \sum_{j=1}^{N}d_je_j(x)=0 $ for all $x\in \supp~\varphi$ which implies that $d_j=0$ for  $j=1,2,\cdots,N$.
Thus the sequence $\{{\cal T}^j c\}_{j=1}^\infty$ is linearly independent, as desired.
\ep

Corollary \ref{60829a} leads to a general result about arbitrary
group representations and the operators generated by cyclic subgroups indexed by $\mz$:

\bc \label{61006b} Let $G$ denote a locally compact group, and $\pi$ a group representation of $G$ on a Hilbert space $\h.$
Given any $x_0\in G$ and any $f_0\in \h,$  and assume that the
family $\{\pi(x_0^k)f_0\}_{k\in \mz}= \{\pi(x_0)^k f_0\}_{k\in \mz}$ is a frame sequence. Then either the
family is a Riesz sequence, or it has infinite excess. \ec
The result in Corollary \ref{61006b} is known in certain special cases, e.g., for the case of a shift-invariant
system considered in Example \ref{motivation example}.

Note that the opposite implication in Corollary \ref{60829a} does not hold; that is,
the operator $T$ is not necessarily bounded even if
$\ftz=\Tnz$ is an overcomplete frame and $\text{dim}(N_U)=\infty$. This is demonstrated by the
following example.

\bex \label{60829b} A collection of functions in $\ltr$ of the form $\mts$ for some $a,b>0$ and some
$g\in \ltr$ is called a {\it Gabor system}.  It is known that if $g\neq 0,$ then
the Gabor system $\mts$ is automatically linearly independent, see
\cite{Linn,Heil3}; thus it can be represented on the form
$\Tnz.$
Now, consider the Gabor frame  $\{E_{m/3}T_{n}\chi_{[0,1]}\}_{m,n\in\mz},$ which is
the union of the three orthonormal bases
$\{E_{k/3}E_{m}T_{n}\chi_{[0,1]}\}_{m,n\in\mz}, \, k=0,1,2.$ The Gabor frame
$\{E_{m/3}T_{n}\chi_{[0,1]}\}_{m,n\in\mz}$
is linearly independent and has infinite excess; in particular $\text{dim}(N_U)=\infty$. Re-order the frame as $\ftz$ in such a way that the elements $\{f_{2k+1}\}_{k\in\mz}$ corresponds to the orthonormal basis $\{E_{m}T_{n}\chi_{[0,1]}\}_{m,n\in\mz}.$ By construction, the elements $\{f_{2k}\}_{k\in\mz}$ now forms an overcomplete frame. By Proposition \ref{60811d}, there is an operator $ T:\Span\ftz\rightarrow \Span\ftz $ such  that $ \ftz = \Tnz $. Since the subsequence $ \{f_{2k}\}_{k\in\mz} $ is an overcomplete frame, there is a non-zero sequence $\{ c_{2k} \}_{k\in\mz}\in \ltz $ such that $ \sum_{k\in\mz}c_{2k}f_{2k}=0 $. Defining $c_k=0 $ for $ k\in 2\mz + 1 $, we have $\sum_{k\in\mz}c_{k}f_{k}=\sum_{k\in\mz}c_{2k}f_{2k}=0$.
On the other hand, since   $ \{f_{2k+1}\}_{k\in\mz} $   is a Riesz basis and $ \ctz $ is non-zero, $ \sum_{k\in\mz}c_{k}f_{k+1}=\sum_{k\in\mz}c_{2k}f_{2k+1} \neq 0 $.  This shows that $N_U$ is not invariant under right-shifts; thus, $T$ is unbounded by Theorem \ref{boundedness}.
\ep \enx

If $\ftz$ is a Riesz basis on the form $=\Tnz,$ then the extension
 of the bounded operator $T: \Span \ftz \to \Span \ftz$ to $\h$ is injective. On the other hand, if a given frame
has the form $\ftz=\Tnz$ for a bounded and injective operator on  $\h$, we can not conclude  that $\ftz$ is a Riesz basis:

\bex  Using the characterization in Example \ref{motivation example}, it is easy to
construct an overcomplete frame sequence $ \{T_k\varphi\}_{k\in\mz}=\{(T_1)^k\varphi\}_{k\in\mz}$
in $\ltr; $ in other words, letting
$\h:= \span \{T_k\varphi\}_{k\in\mz}$ the sequence  $ \{T_k\varphi\}_{k\in\mz}$ is an overcomplete
frame for $\h.$ Clearly $T_1$ is a bounded and injective operator on $\h,$ but by construction
 $ \{T_k\varphi\}_{k\in\mz}$ is not a Riesz basis for $\h$.
\ep \enx

			\bc Consider a tight frame having a representation $\ftz=\Tnz$
 for some invertible operator $T\in B(\h).$ Then $T$ is an isometry.
			\ec\bp Since the frame bounds are $A=B$, using Theorem \ref{boundedness}, we have $\|T\|=\|T^{-1}\|=1$. Therefore
			$ \|f\|=\|T^{-1}Tf\|\leq\|Tf\|\leq\|f\|, $ which implies that $T$ is in isometry. \ep

 \section{Duality} \label{61003b}

In this section we will analyze certain aspects of the duality theory for a  frame having the form $\ftz=\Tnz$
for some bounded linear and invertible operator $T: \Span \ftz\to \Span \ftz.$ In particular
we will identify a class of dual frames (including the canonical dual frame)
that is also given by iteration of a bounded operator. On the other hand, we also give
an example of a frame for which not all dual frames have this form.

In the entire section we denote
the synthesis operator by $U$; then the frame operator is $S=UU^*.$
We first prove
that the synthesis operator $U$  is an intertwining operator for the right-shift operator ${\cal T}$
on $\ltz$ and the operator $T$, as well as an immediate consequence for the frame operator.
Let $c_{00} \subset \ltz$ denote the subspace consisting of finite sequences.

	\bl Consider a  Bessel sequence having the form $\ftz=\Tnz$
for a linear operator $T: \Span \ftz\to \Span \ftz.$ Then $TU= U{\cal T}$ on  $c_{00}.$
%\bes TU= U{\cal T} \, \, \mbox{on} \, \, c_{00}.\ens
Assuming that $T$ has an extension to a
bounded operator $\widetilde{T}:\h \to \h,$ the following hold:
\bei \item[(i)]  $\widetilde{T}U= U{\cal T}$ on $\ltz$.
\item[(ii)]\label{lem} If $\Tnz$ is a frame and $\widetilde{T}$ is invertible, then
		$\widetilde{T}S=S(\widetilde{T}^*)^{-1};$ in particular,
$ S\wt = \wt S $ if and only if $\wt$ is unitary.
\eni
	\el
	\bp For $\ctz\in c_{00}$,  there is an $N\in\mn$ such that $c_k=0$ for $|k|\geq N$. Therefore
	\bee \label{intertwining} TU\ctz &=& T\sum_{k=-N}^N c_kf_k=\sum_{k=-N}^N c_{k}f_{k+1}=\sum_{k=-N+1}^{N+1} c_{k-1}f_k \nonumber\\
	&=&U\{c_{k-1}\}_{k\in\mz}=U{\cal T}\ctz.\ene
	In the case that $\widetilde{T}$ is bounded, the equality \eqref{intertwining} holds on $\ltz$ because $c_{00}$ is dense in $\ltz$; this proves (i).  For the proof of (ii), using \eqref{intertwining}
and  that  $S=UU^*$,
		\[\widetilde{T}S\widetilde{T}^*=\widetilde{T}UU^*\widetilde{T}^*=\widetilde{T}U(\widetilde{T}U)^*
=U{\cal T}(U{\cal T})^*=U{\cal T}{\cal T}^*U^*=UU^*=S.\]
		Therefore $\widetilde{T}S=S(\widetilde{T}^*)^{-1}$, as desired.
	\ep
% \begin{rem}\label{besselintertwining}
% 	 For any Bessel sequence $\ftz=\Tnz$, Lemma \ref{intertwining} (i)  holds on a dense subset of $ \ltz $ even if $T$ is not bounded on $\h$ . Denote the subspace of all finite sequence by $ c_{00} $. Then
% \end{rem}

For a frame $\ftz= \Tnz,$ the operator $S^{-1}TS$ is
invertible considered  from $\Span \{S^{-1}f_k\}_{k\in\mz}$ into itself, and  the canonical dual frame is $\{S^{-1}f_k\}_{k\in\mz}=\{(S^{-1}TS)^kS^{-1}f_0\}_{k\in\mz}.$ This was already observed  in the finite-dimensional setting in \cite{AK}. In the case where  $ T $ has an extension to a bounded and invertible
operator on $\h$ (see the appropriate conditions in Theorem \ref{boundedness}), we will now derive an alternative description of the canonical dual frame, directly in terms of the operator $ T$ and its adjoint. Since the rest of the results in the current section will use the
same assumptions on the operator $T,$ we will drop the distinction between the operator
$T$ and $\widetilde{T},$ and simply denote the operator by $T.$

	\begin{prop}\label{canonicaldual}
		Consider a frame $\ftz=\Tnz$, where $T\in B(\h)$
is invertible. Let  $\widetilde{f_0}=S^{-1}f_0$. Then  $\{S^{-1}f_k\}_{k\in\mz}=\{(T^*)^{-k}\widetilde{f_0}\}_{k\in\mz}$.
	\end{prop}
	\bp Lemma \ref{lem} (ii) implies that $TS=S(T^*)^{-1}. $ Thus $S^{-1}T=(T^*)^{-1}S^{-1}$ and therefore $S^{-1}T^k=(T^*)^{-k}S^{-1}$ for $k\in\mn$. We also have  that
		 $S^{-1}T^{-1}=T^*S^{-1}$ and thus  $S^{-1}T^{-k}=(T^*)^kS^{-1}$ for $k\in\mn$. It follows that
\bes \{S^{-1}f_k\}_{k\in\mz}=\{S^{-1}T^kf_0\}_{k\in\mz}=\{(T^*)^{-k} S^{-1}f_0\}_{k\in\mz}=\{(T^*)^{-k} \widetilde{f_0}\}_{k\in\mz},\ens as desired.
	\ep

%	\bl Assume that $\Tnz$ is a frame with synthesis operator $U$. The following hold:
%		\bei
%		\item[(1)]\label{1lem} For all $h_0\in\Span\Tnz$, the sequence $\{T^kh_0\}_{k\in\mz}$ is frame.
%		\item[(2)]\label{2lem}  $\{g_k\}_{k\in\mz}$ is a dual frame if and only if for a bounded operator $W:\ltz\rightarrow\h$, we have  $g_j=(S^{-1}U+ W(I-U^*S^{-1}U))\delta_j$.
%		\item[(3)]\label{4lem} The sequence $\{(S^{-1}TS)^k \widetilde{f_0}\}_{k\in\mz}$ is frame for $\widetilde{f_0}=S^{-1}f_0$.
%		\eni
%		\el

Since the translation operators on $\ltr$ are unitary,
Proposition \ref{canonicaldual} generalizes the well-known result that the canonical dual
of a shift-invariant frame $\{T_k\varphi\}_{k\in\mz}$ in $\ltr$ has the form
$\{T_k\widetilde{\varphi}\}_{k\in\mz}$ for some $\widetilde{\varphi}\in\ltr$

It is important to notice that Proposition \ref{canonicaldual} only shows that the canonical dual frame has the form of an iterated system. Indeed,
the next example exhibits a frame satisfying the conditions
in Proposition \ref{canonicaldual}   and having a dual frame that is not representable
by an operator:

\bex \label{61003f}Let us return to Example \ref{motivation example} and consider
an overcomplete frame sequence $\{T_k\varphi\}_{k\in\mz}$ in $\ltr. $ Then there exists an
element $T_{k^\prime}\varphi, \, k^\prime \in \mz,$ that can  be removed from
the frame sequence, leaving a frame sequence for the same space; due to the special
structure of the frame we can even take $k^\prime =0.$ Letting $\{g_k\}_{k=-\infty}^{-1}
\cup \{g_k\}_{k=1}^\infty$ denote a dual frame for the resulting frame sequence
$ \{T_k\varphi\}_{k=-\infty}^{-1} \cup \{T_k\varphi\}_{k=1}^\infty$ this implies that
the frame $\{T_k\varphi\}_{k\in\mz} $ has the non-canonical dual
 $\{g_k\}_{k=-\infty}^{-1} \cup \{0\}
\cup \{g_k\}_{k=1}^\infty;$ this family is clearly linearly dependent. Hence, by
Theorem \ref{60811d} the system is
not representable by an operator.
\ep \enx

We will now show that despite the obstruction in Example \ref{61003f} we can
 actually characterize the class  of dual frames that arise through
 iterated actions of a bounded operator. We first show that the
 only candidate for this operator indeed is the operator  $(T^*)^{-1}$ arising
 in Proposition \ref{canonicaldual}. In particular, this shows that for a frame $\ftz=\Tnz$ given in terms
 of a unitary operator $T,$ the dual frames having the form of an iterated operator system must be
 generated by the same operator.
 \bl\label{61005a} Consider a frame $\ftz=\Tnz$, where $T\in B(\h)$
is invertible. Assume that  $\gtz=\Vnz$ is a dual frame and that $V$ is bounded. Then $V=(T^*)^{-1}.$
 \el
 \bp For any $f\in B(\h)$, two applications of the frame decomposition yield that
\bes
 f &=& \sum_{k\in\mz}\langle f,T^kf_0 \rangle V^kg_0=V\sum_{k\in\mz}\langle f,T^{k}f_0 \rangle V^{k-1}g_0\\
 &=& V\sum_{k\in\mz}\langle T^*f,T^{k-1}f_0 \rangle V^{k-1}g_0
 =V\sum_{k\in\mz}\langle T^*f,T^{k}f_0 \rangle V^{k}g_0=VT^*f.
\ens
Therefore $VT^*=I$. Since $T$ is invertible it follows that $V=(T^*)^{-1}$.
\ep

We will now give the full characterization of dual frames of $\ftz=\Tnz$ that
are given in terms of iterations of a bounded operator.

			\begin{thm} \label{60829g}
Consider a frame $\ftz=\Tnz,$ where $T\in B(\h)$ is invertible.
Then the dual frames given as iterates of a bounded operator are precisely the families
of the form $\{(T^*)^{-k} g_0\}_{k\in\mz}$  for which
\bee \label{61001a} g_0=S^{-1}f_0+h_0-\sum_{j\in\mz}\langle S^{-1}f_0,T^jf_0\rangle (T^*)^{-j}h_0\ene for some $h_0\in \h$ such that $\{(T^*)^{-k}h_0\}_{k\in \mz}$
is a Bessel sequence.
In particular, this condition is satisfied when $h_0$ is taken
from the dense subspace $\mbox{span} \{(T^*)^{-k}\widetilde{f_0}\}_{k\in \mz}.$
			\end{thm}
			\bp First, note that by Lemma \ref{61005a} we know
that the only operator that might be applicable in
the desired representation of the dual frame is $(T^*)^{-1}.$ Now,   assume that
$\{(T^*)^{-k} \varphi\}_{k\in\mz}$ is a dual frame of $\Tnz$
for some $\varphi\in \h.$ Then $\{(T^*)^{-k} \varphi\}_{k\in\mz}$
is a Bessel sequence, and  taking $h_0:=\varphi$ in \eqref{61001a} yields
that $g_0=\varphi$. On the other hand,
assume that $h_0\in \h$ is chosen such that the sequence $\{(T^*)^{-k} h_0\}_{k\in\mz}$ is a Bessel sequence, and choose $g_0$ as in \eqref{61001a}.  Denote the synthesis operator of $\{(T^*)^{-k} h_0\}_{k\in\mz}$ by $W$. Letting
$\{\delta_j\}_{j\in \mz}$ denote the canonical orthonormal basis for
$\ltz$ and $ V:=S^{-1}U+ W(I-U^*S^{-1}U),$  it follows
 from \cite{Li} (alternatively,
see Lemma 6.3.5 and Lemma 6.3.6 in \cite{CB}) that
the sequence $\{V\delta_j\}_{j\in\mz}$ is a dual frame of $\Tnz$.
 Furthermore, by direct calculation,
\bes V \delta_0=  \left(S^{-1}U+ W(I-U^*S^{-1}U)\right)\delta_0=       g_0.\ens
We first show that this
frame indeed has the form $\{(T^*)^{-k} g_0\}_{k\in\mz}.$  We will now show that $V$ is an intertwining operator between ${\cal T}$ and $(T^*)^{-1}$.
Applying Lemma \ref{lem}  (i) on $\{(T^*)^{-k} h_0\}_{k\in\mz}$, we know that $W{\cal T}=(T^*)^{-1}W$ and $TU=U{\cal T}$.  Also since $T$ is bounded and invertible, Lemma \ref{intertwining} (ii) shows that  $S^{-1}T=(T^*)^{-1}S^{-1}.$  Hence we get
			\bes
			V{\cal T}
			&=&(S^{-1}U+ W(I-U^*S^{-1}U)){\cal T}\\
			& = & S^{-1}U{\cal T}+ W{\cal T}-WU^*S^{-1}U{\cal T}\\
			&=& S^{-1}TU+ (T^*)^{-1}W-WU^*S^{-1}TU\\
			&=&(T^*)^{-1}S^{-1}U+ (T^*)^{-1}W-WU^*(T^*)^{-1}S^{-1}U.
			\ens
			Similarly to the proof of Lemma \ref{lem}, $T^{-1}U=U{\cal T}^{-1}$. Therefore $U^*(T^*)^{-1}=(T^{-1}U)^*=(U{\cal T}^{-1})^*={\cal T}U^*$; thus,
			
			\[V{\cal T}=(T^*)^{-1}(S^{-1}U+ W-WU^*S^{-1}U)=(T^*)^{-1}V.\] This implies that $V\delta_j=V{\cal T}^j\delta_0 = (T^*)^{-j}V\delta_0
			= (T^*)^{-j} g_0,$ as desired.

Finally, we note that if $h_0\in\Span \{(T^*)^{-k} \widetilde{f_0}\}_{k\in\mz}$, the sequence $\{(T^*)^{-k} h_0\}_{k\in\mz}$ is a finite  sum of frame sequences and hence  a Bessel sequence.\ep

In order to apply Proposition \ref{canonicaldual} and Theorem \ref{60829g} we must
calculate the adjoint of the operator $T$ arising in the representation $\ftz=\Tnz.$
In general this can only be done with specific knowledge of the operator $T$ at hand.
An
additional condition on the frame $\ftz$ implies that the operator $T$ is unitary,
and allows us to find it explicitly in terms of $\ftz;$ the result generalizes
the observations for shift-invariant
systems in
Example \ref{motivation example}, and also applies to some of the other systems obtained via group
representations in Example \ref{61006a}.

 	\begin{prop}\label{Z}Consider a  frame having the form $\ftz=\Tnz$
for some  operator $T\in B(\h).$
 		Assume that for a function $ \theta:\mz\mapsto \R $, we have
 		$\langle f_j,f_k \rangle=\theta(j-k),\quad  j,k\in \mz.$ Then
 	$ \widetilde{T}^* \sum_{k\in\mz} c_k f_k = \sum_{k\in\mz} c_k f_{k-1}$ for all $\{ c_k\}_{k\in\mz}\in\ell^2(\mz). $ In particular, $\widetilde{T}$ is unitary.
 	\end{prop}
 	\bp
 		Consider  arbitrary $j,k\in\mz$. Then
 		\bes
 		\langle \widetilde{T}f_j,f_k\rangle =\langle f_{j+1},f_{k}\rangle
 		=\theta(j+1-k)=\langle f_j,f_{k-1}\rangle= \langle f_j,\widetilde{T}^*f_k\rangle.
 		\ens
 		It follows that $\widetilde{T}^*f_k=f_{k-1}$. Therefore
 $\widetilde{T}\widetilde{T}^*=\widetilde{T}^*\widetilde{T}=I,$ i.e., $\widetilde{T}$ is unitary and  $$ \widetilde{T}^* \sum_{k\in\mz} c_k f_k = \sum_{k\in\mz} c_k f_{k-1}$$ for all $\ctz\in \ltz.$
 	\ep

\section{Stability of the representation $\ftz=\Tnz$} \label{61003c}
For applications of frames it is important that key properties are kept under
perturbations. We will now
state a perturbation condition that preserves the existence of
a representation $\ftz=\Tnz$.  The condition was first used in connection with frames in
the paper \cite{OC}.

\begin{prop} \label{61001b}
	Assume that $\ftz=\Tnz$ is a frame for $\h$ and let $\gtz$ be a sequence in $\h$. Assume that  there exist constants $\lambda_1,\lambda_2\in [0, 1[$ such that
	\begin{equation}\label{lambda1,2}
	\nl \sum c_k(f_k-g_k)   \nr\leq \lambda_1 \nl \sum c_kf_k \nr +\lambda_2 \nl \sum c_kg_k \nr	
	\end{equation}
for all finite sequences $\{c_k\}.$ Then $\gtz$ is  a frame for $\h$;
furthermore $\gtz$ can be represented as $\gtz=\Vnz$ for a linear operator
\bes V:\Span\gtz\rightarrow\Span\gtz.\ens
If $T$ is bounded, then $V$ is also bounded.
\end{prop}
\bp
By Theorem 2 in \cite{OC} the perturbation condition implies that $\gtz$ is a frame. Also, since $\max(\lambda_1,\lambda_2)<1$, it follows from \eqref{lambda1,2}  that
\begin{equation}\label{ifandonlyif}
	\sum_{k\in\mz}c_kf_k=0\quad\Leftrightarrow \quad \sum_{k\in\mz}c_kg_k=0,\quad
\forall \ctz\in\ltz.
\end{equation}
Since $\ftz$ is linearly independent, \eqref{ifandonlyif} implies that the sequence  $\gtz$  also is linear independent. Therefore by Proposition \ref{60811d}, there is a linear operator $V:\Span\gtz\rightarrow\Span\gtz$ such that
$\gtz=\Vnz$. Now assume that the operator $T$ is bounded; We want to show that then $V$ is also bounded. Let $W:\ltz\rightarrow\h$ be the synthesis operator for $\gtz$,
and consider some $\ctz\in N_W$. Then by \eqref{ifandonlyif},  $\ctz \in N_U$, where $U$ is the synthesis operator for $\ftz$. Since $T$ is bounded, Theorem \ref{boundedness}
implies that
$N_U$ is invariant under right-shifts, i.e., $\sum_{k\in\mz}c_{k-1}f_k=0$. Using again \eqref{ifandonlyif}, we conclude that $\sum_{k\in\mz}c_{k-1}g_k=0$, which shows that ${\cal T}\ctz\in N_W$. Applying Theorem \ref{boundedness} again shows that  $V$ is bounded.
\ep

Note that \eqref{lambda1,2} is a special case of the  perturbation condition
\begin{equation}\label{lambdamu}
\nl\sum c_k(f_k-g_k)   \nr\leq \lambda_1 \nl \sum c_kf_k \nr+ \lambda_2 \nl \sum c_kg_k \nr + \mu (\sum|c_k|^2)^{1/2},\end{equation}
appearing in \cite{OC}. If $\ftz$ is a frame
for $\h$ with lower bound $A,$ $\gtz\subset \h,$
and \eqref{lambdamu} holds for all finite sequences $\{c_k\}$ and some parameters $\lambda_1, \lambda_2, \mu\ge 0$
such that $\mbox{max}\left(\lambda_2, \lambda_1+ \frac{\mu}{\sqrt{A}}\right)<1,$ then
by \cite{OC} also $\gtz$ is a frame for $\h.$
This perturbation condition has been used in many different contexts in frame theory, typically
for the case $\mu>0.$ However,  the case $\mu>0$ turns out to be problematic if we want the perturbation $\gtz$
of a frame $\ftz= \Tnz$
to be represented on the form $\gtz= \{W^kg_0\}_{k\in \mz}.$
The first obstacle is that if $\mu>0,$ the perturbation condition \eqref{lambdamu}     does not preserve the
property  of being representable by an operator:

\bex
Consider an orthonormal basis $\{e_k\}_{k\in\mz}$ for a Hilbert space $\h$.
Then the family
$\{f_k\}_{k\in I}:= \{e_k\}_{k\in\mz} \cup \{ \alpha \sum_{j=1}^\infty \frac1{2^j}\, e_j\}$
is a linearly independent frame for any choice of $\alpha> 0,$ with lower frame
bound $A=1.$
For $\alpha <1,$ the family
$\{g_k\}_{k\in I}:= \{e_k\}_{k\in\mz} \cup \{ 0\}$
is a perturbation of $\{f_k\}_{k\in I}$ in the sense of \eqref{lambdamu}, with
$\lambda_1=\lambda_2=0$ and $\mu =\alpha.$ However, regardless how small we choose
$\alpha,$ the family $\{g_k\}_{k\in I}$ is not linearly independent.
Hence, by Proposition \ref{60811d} the sequence $\{g_k\}_{k\in I}$ can not be
represented on the form $\{W^k\varphi\}_{k\in \mz}.$ \ep
\enx

The following example shows that even if we assume that the perturbation $\gtz$
of a frame $\ftz= \Tnz$ is linearly independent
(and hence representable on the form $\gtz= \{V^k g_0\}_{k\in \mz}$), the condition
\eqref{lambdamu} does not imply that $V$ is bounded if $T$ is bounded.

\bex \label{61017a} Let us first explain the idea of the construction in the setting of a general Hilbert space
$\h$.  Assume that $\ftz=\Tnz$ is an overcomplete frame for $\h,$
with lower bound $A,$ and that the operator $T$ is bounded. We further assume that
\bei
\item[(a)] The sequence $\{f_k\}_{k\in \mz \setminus \{-1,0\}}$ is complete in $\h.$
\eni
We will then search for some $g_0\in \h$ such that the sequence
\bes \gtz:= \{f_k\}_{k=-\infty}^{-1} \cup \{g_0\} \cup \{f_k \}_{k=1}^\infty\ens satisfy the following requirements:
\bei
\item[(b)] The condition \eqref{lambdamu} is satisfied with
$\mbox{max}\left(\lambda_2, \lambda_1+ \frac{\mu}{\sqrt{A}}\right)<1;$
\item[(c)] $\gtz$ is linearly independent.
\eni  We will now explain how this setup leads to the desired conclusion; after that we provide a
concrete construction satisfying all the requirements.

First, the condition (b) implies that $\gtk$ is a frame for $\h;$ by (c) it has the form
$\gtz= \{W^kg_0\}_{k\in \mz}$ for some operator  \bes W: \mbox{span} \gtz \to \mbox{span} \gtz.\ens
By the definition of the sequence $\gtz$ it follows that

\bee \label{61010af} \left\{\begin{array}{lll}
 	 	Wf_k &= f_{k+1}, \hspace{.3cm} k=-2,-3, \dots\\
 	 	Wf_{-1} &= g_0,\\
 Wg_{0} &=f_1, \\
 Wf_{k} &= f_{k+1},\hspace{.3cm} k=1,2, \dots
 	 \end{array}\right.\ene
We note that the operators $T$ and $W$ act in an identical way on the vectors
 $\{f_k\}_{k\in \mz \setminus \{-1,0\}};$ thus, if $W$ is bounded it follows by (a) that $W=T.$
 But then \eqref{61010af} implies that $g_0=Wf_{-1}= Tf_{-1}=f_0,$ i.e., that
 $\ftz=\gtz.$ In other words: for a perturbation satisfying the stated conditions, the operator
 $W$ in the representation $\gtz= \{W^kg_0\}_{k\in \mz}$ will not be bounded when $g_0 \neq f_0.$

We now proceed to a concrete construction satisfying (a)--(c).
In order to do so, we return to the shift-invariant systems
considered in Example \ref{motivation example}.
First, it is well-known that the function $\mbox{sinc}(x):= \frac{\sin (\pi x)}{\pi x}$
generates an orthonormal basis $\{T_k \mbox{sinc}\}_{k\in \mz}$ for the Paley-Wiener space
$$\h:= \{f\in \ltr \, \big| \, \mbox{supp}\, \widehat{f} \subseteq [-1/2, 1/2]\}.$$ It follows that
the oversampled family $\ftz:=\{T_{k/3} \mbox{sinc}\}_{k\in \mz}=\{T_{1/3}^k \mbox{sinc}\}_{k\in \mz}$ can be considered as a union of three
orthonormal bases, and hence form a tight frame for $\h;$ we note that by the carefully chosen
oversampling, the condition (a) is satisfied. The operator $T:= T_{1/3}$ is clearly bounded.

Now, consider a constant $c\ge 0$ and let $g_0:= T_{c}f_0;$ then $g_0\in \h,$ and   for any
finite scalar sequence $\{c_k\}$ we have
\bes \nl \sum c_k (f_k -g_k) \nr = \nl c_0 (f_0-T_cf_0)\nr
\le \nl f_0-T_cf_0 \nr   (\sum|c_k|^2)^{1/2}.\ens
By continuity of the translation operator there exists some $\delta>0$ such that
$\nl f_0-T_cf_0 \nr < \sqrt{A}$ whenever $c\in [0, \delta[;$ it now follows from the
perturbation condition that $\gtz$ is a frame for $\h$ for $c$ belonging to this range, i.e., the
condition (b) is satisfied.
Furthermore, for $c< 1/3$ all the translation parameters appearing in the sequence $\gtz$ are
pairwise different; thus $\gtz$ is linearly independent and condition (c) is fulfilled.
\ep \enx

Most of the concrete applications of perturbation results in frame theory deals with
the special case of the condition \eqref{lambdamu} corresponding to $\lambda_1=\lambda_2=0.$
Even in this case, Example \ref{61017a} shows that the perturbation condition does not
preserve boundedness of the representing operator. In applications where stability is an important issue,
 one can alternatively represent a frame using iterated operator systems based on a finite collection of
operators instead of a singleton. Consider a frame $\ftz$ which is
norm-bounded below.
It is proved in \cite{olemmaarzieh} that then there is a finite collection of
vectors from $\ftz,$ to be called $\varphi_1, \dots, \varphi_J,$ and
corresponding bounded operators $T_j: \h \to \h,$ such that
$\{T_j^n \varphi_j\}_{n\in \mz}$ is a Riesz sequence, and
$ \label{50614a} \ftz = \cup_{j=1}^J \{T_j^n \varphi_j\}_{n\in \mz}.$ The proof uses
the Feichtinger theorem (which was
a conjecture for several years and finally got confirmed in \cite{MSS}).
We will now show that the stated representation  is
stable with respect to the central perturbation condition \eqref{lambdamu} with $\lambda_1=\lambda_2=0.$
	
\bt \label{61017b} Assume that the frame $\ftz$ is norm-bounded below, and consider
a representation on the form $\cup_{j=1}^J \{T_j^n \varphi_j\}_{n\in \mz},$ where the operators
$T_j$ are bounded and $\{T_j^n \varphi_j\}_{n\in \mz}$ is a Riesz sequence.
 Let $A$ denote a common lower frame bound
for the frame $\ftk$ and all the Riesz sequences $\{T_j^n \varphi_j\}_{n\in \mz}, \, \, j=1, \dots, J.$ Let $\gtz$
be a sequence in $\h$ such that for some $\mu < \sqrt{A},$
\begin{equation}\label{mu}
\nl\sum c_k(f_k-g_k)   \nr\leq  \mu (\sum|c_k|^2)^{1/2},\end{equation}
for all finite scalar sequences $\{c_k\}.$ Then there is a finite collection of
vectors $\phi_1, \dots, \phi_J$ from $\gtz$ and
corresponding bounded operators $W_j: \h \to \h,$ such that
$\{W_j^n \phi_j\}_{n\in \mz}$ is a Riesz sequence, and
\bee \label{50614ag} \gtz = \bigcup_{j=1}^J \{W_j^n \phi_j\}_{n\in \mz}.\ene
\et

\bp The perturbation condition \eqref{mu}  is a special case of \eqref{lambdamu}; thus
the family $\gtk$ is a frame for $\h.$ Furthermore, partitioning $\gtz$
according to the splitting $\ftz=\cup_{j=1}^J \{T_j^n \varphi_j\}_{n\in \mz},$ i.e., writing
$\gtz = \bigcup_{j=1}^J \{ g_j^{(n)}\}_{n\in \mz},$
it follows from \eqref{mu} that for any fixed $j\in \{1, \dots, J\}$ and any finite scalar sequence $\{c_n\}_{n\in \mz},$
 \bes \nl \sum c_n(T_j^n \varphi_j-g_j^{(n)} )\nr\leq  \mu (\sum|c_n|^2)^{1/2}. \ens
 Therefore, for each fixed $j\in \{1, \dots, J\}$ the sequence $\{ g_j^{(n)}\}_{n\in \mz}$ is
 a Riesz sequence, and hence representable on the form $\{W_j^n \phi_j\}_{n\in \mz}$
 for some bounded operator $W_j$ and some $\phi_j\in \h.$

Note that an alternative way of proving the result would be to show directly that
\eqref{mu} implies that $\gtz$ is norm-bounded below and then refer to the stated result in
\cite{olemmaarzieh}.
However, this argument would not yield that the splitting of the indexing of the frame $\ftz$ is preserved for the perturbed family $\gtz,$ as in \eqref{50614ag}. \ep

\section*{Appendix: auxiliary examples}
We will close the paper with a few operator-theoretical considerations, to which we have referred
throughout the paper.

\noindent {\bf 1)} Instead of representing a frame on the form
$\Tnz,$ one could also consider
representations on the form $\{T^nf_0\}_{n=0}^\infty;$ this indexing occur,
e.g., in dynamical sampling \cite{A1,A3}. The chosen indexing actually
has a serious influence on the properties of the operator $T.$ For
example, there exist frames  $\ftz= \Tnz$ where $T$ is a unitary operator (Example \ref{motivation example}); but if we reindex the frame
as $\{f_k\}_{k=0}^\infty$ it can not be represented on the form $\{U^nf_0\}_{n=0}^\infty$ for a unitary operator $U,$ see \cite{A3}.

Let us demonstrate the sensibility to the indexing by one more case.
First,
it is well-known that the function $\mbox{sinc}(x):= \frac{\sin (\pi x)}{\pi x}$
generates an orthonormal basis $\{T_k \mbox{sinc}\}_{k\in \mz}$ for the Paley-Wiener space
$$\h:= \{f\in \ltr \, \big| \, \mbox{supp}\, \widehat{f} \subseteq [-1/2, 1/2]\}.$$
It follows that
the oversampled family $\{T_{k/2} \mbox{sinc}\}_{k\in \mz}=\{T_{1/2}^k \mbox{sinc}\}_{k\in \mz}$
is a tight frame for $\h;$ the representing
operator $T_{1/2}$ is clearly bounded.

On the other hand, the following lemma shows that considering $\{T_{k/2} \mbox{sinc}\}_{k\in \mz}$
as a union of the two orthonormal bases $\{T_{k} \mbox{sinc}\}_{k\in \mz}$ and
$\{T_{1/2}T_{k} \mbox{sinc}\}_{k\in \mz},$ re-indexing
in a natural fashion as $\{T^nf_0\}_{n=0}^\infty$ always leads to an unbounded operator $T.$

\bl Consider two orthonormal bases $\ftk$ and $\etk$ for a Hilbert space $\h$ and assume that the set
$\ftk \cup \etk$ is linearly independent.
Then  a linear operator $T: \mbox{span}\left(\ftk \cup \etk\right) \to \h$ such that
\bee \label{61124a} \{\varphi_k\}_{k=0}^\infty:= \{f_1,e_1, f_2, e_2, \dots\}=\{T^k \varphi_0\}_{k=0}^\infty,\ene
is necessarily  unbounded. \el

\bp The ordering in \eqref{61124a} implies that
$Tf_k=e_k$ for all $k\in \mn;$ thus, if the operator $T$  is bounded, it has a unique
extension to a bounded linear operator
$\widetilde{T}: \h \to\h$, given by
$\widetilde{T} \sum_{k=1}^\infty c_kf_k = \sum_{k=1}^\infty c_ke_k, \, \ctk \in \ltn.$
Clearly $\widetilde{T}$ is a surjective mapping. On the other hand, \eqref{61124a} also
implies that $\widetilde{T}e_k= f_{k+1};$ since $\etk$ is an orthonormal basis for $\h$
this implies that the range of $T$ equals the space $\span \{f_k\}_{k=2}^\infty,$
which excludes that $\widetilde{T} $ is surjective.
This contradiction shows that the operator $T$ can not be bounded. \ep

\noindent {\bf 2)} Let $V$ denote a dense subspace of a Hilbert space $\h,$
and consider a bounded and bijective operator $T: V\to V.$ Then $T$ has a unique  extension
to a bounded operator $\widetilde{T}: \h\to \h.$  The following example demonstrates that
the
extension $\widetilde{T}$ might no longer be injective.
\bex \label{60828h} Let $\etk$ denote an orthonormal basis for a separable Hilbert space,
and consider the sequence $\ftk:= \{e_1\} \cup \{e_{k-1}+\frac1{k}\, e_k\}_{k=2}^\infty.$
Then $\ftk$ is a frame (see \cite{CB}), and it is easy to see that  the elements are linearly independent.
Let $V:=\Span \ftk,$ and consider the operator
$T: V\to V,\, Tf:= \sum_{k=1}^\infty \la f, e_k\ra f_k.$
Then $T$ is linear and bounded, its extension to a bounded operator on $\h$ is given by the
same expression, i.e.,
$\widetilde{T}: \h \to \h, \, \widetilde{T}f= \sum_{k=1}^\infty \la f, e_k\ra f_k.$
The operator $T:V \to V$ is bijective. Injectivity follows from the fact that for
$f\in V,$ the sequence $\{\la f, e_k\ra\}_{k=1}^\infty$ is finite; hence, due to the
linear independence of $\ftk$ we can only have $Tf=0$ if $\la f, e_k\ra=0$ for all $k\in \mn,$
which implies that $f=0.$  So show that $T$ is surjective, let $g\in V.$ Then
$g=\sum_{k=1}^N c_kf_k$ for some $N\in \mn$ and some $c_k\in \mc.$ Then $f:= \sum_{k=1}^N c_ke_k\in V$
and $Tf=g.$

However, the bounded extension $\widetilde{T}:\h \to \h$ is not injective. Indeed, since
$\ftk$ is an overcomplete frame, there exist coefficients $\ctk\in \ltn \setminus \{0\}$ such
that $\sum_{k=1}^\infty c_kf_k=0;$ taking $f:= \sum_{k=1}^\infty c_ke_k\neq 0$ we
have $Tf=0.$\ep \enx

\begin{tabbing}
text-text-text-text-text-text-text-text-text-text \= text \kill \\
Ole Christensen \> Marzieh Hasannasab \\
Technical University of Denmark \> Technical University of Denmark  \\
DTU Compute \> DTU Compute \\
Building 303, 2800 Lyngby \> Building 303, 2800 Lyngby \\
Denmark \> Denmark \\
Email: ochr@dtu.dk \> mhas@dtu.dk
\end{tabbing}

\end{document}